\documentclass{amsart}

\usepackage{amsaddr}
\usepackage[numbers]{natbib}
\usepackage{float}
\usepackage{hyperref}
\usepackage[pdftex]{graphicx}
\usepackage{algorithm}

\newtheorem*{prop}{Proposition}
\theoremstyle{remark}\newtheorem*{rmk}{Remark}

\begin{document}

\title{A Note On The Randomized Kaczmarz Method With A Partially Weighted Selection Step}

\author{J\"{u}rgen Gro{\ss}}
\address{Institute for Mathematics and Applied Informatics, University of Hildesheim, Germany}
\email{juergen.gross@uni-hildesheim.de}



\date{}

\begin{abstract} In this note we reconsider two known algorithms which both usually converge faster than the randomized Kaczmarz method introduced by \citet{strohmer2009randomized}, but require the additional computation of all residuals of an iteration at each step. As already indicated in the literature, e.g. \citet{steinerberger2020weighted, jiang2020kaczmarz}, it is shown that the non-randomized version of the two algorithms converges at least as fast as the randomized version, while still requiring computation of all residuals. Based on that observation, a new simple random sample selection scheme has been introduced  by \citet{jiang2020kaczmarz} to reduce the required total of residuals. In the same light we propose an alternative random selection scheme which can easily be included as a `partially weighted selection step' into the classical randomized Kaczmarz algorithm without much ado. Numerical examples show that the randomly determined number of required residuals can be quite moderate.
\end{abstract}

\maketitle

\markboth{J.~Gro{\ss}}{Partially Weighted Randomized Kaczmarz}

\section{Introduction}\label{sec:intro}

Consider the system of linear equations denoted by
\begin{equation}\label{E1}
A x = b
\end{equation}
where $A\in \mathbb{C}^{m\times n}$ for $m\geq n$ is a matrix of rank $n$.

Let $a_{i}$ denotes the $i$-th row of the matrix $A$ considered as a column vector. Starting with an initial guess $x_{0}$, the simple \citet{karczmarz1937angenaherte} method provides an iterative algorithm for the solution to (\ref{E1}) with respect to $x$.
The $(k+1)$-th iteration is given by
\begin{equation}\label{E2}
x_{k+1} = x_{k} + \frac{b_{i} - \langle a_{i}, x_{k}\rangle}{\|a_{i}\|_{2}^2} a_{i}\; ,
\end{equation}
where $i = (k \mod m) +1$, and $\|\cdot \|_{2}$ denotes the Euclidean norm of a vector. The algorithm generates cycles of $m$ iterations by sweeping repeatedly through the rows of the matrix $A$ with the last iteration from one cycle as the starting value for the first iteration of the next cycle. This procedure is also known as the algebraic reconstruction technique (ART), e.g. \citet{gordon1970algebraic}.

The randomized Kaczmarz method by \citet{strohmer2009randomized} improves upon the simple Kaczmarz method by
choosing row $i$ in (\ref{E2}) not in a successive manner but randomly from the set of row numbers $\{1,\ldots , m\}$, where row $i$ is given probability $\|a_{i}\|_{2}^{2}/\|A\|_{F}^2$, and $\|A\|_{F}$ denotes the Frobenius norm of $A$.

Convergence rates of Kaczmarz related methods and procedures for improvements were investigated by a number of authors, see e.g. \cite{deutsch1984rate, deutsch1997rate, galantai2005rate, zouzias2013randomized, elfving2014semi, needell2014paved, gower2015randomized, ma2015convergence, needell2015randomized, liu2016accelerated, hefny2017rows, bai2018greedy, bai2018relaxed, popa2018convergence, gower2019adaptive, zhang2019new, guan2020note, jiang2020kaczmarz, steinerberger2020weighted}.

In the following we will assume that the matrix $A$ has Euclidean row norm \begin{equation}
\|a_{i}\|_{2}=1\end{equation}
for $i=1,\ldots, m$, in which case the matrix $A$ is also called {\em standardized}, see
\citet{needell2014paved}. In that case $\|A\|_{F}^2 = m$ and row $i$ is selected by the randomized Kaczmarz method with  uniform probability $1/m$.

\section{Two Algorithms} \label{sec:algo}

Recently, algorithms had been proposed which randomly select a row for computing $x_{k+1}$ by assigning a probability to row $i$ depending on the $i$-th element
\begin{equation}\label{E3}
r_{k}(i) :=
b_{i} - \langle a_{i}, x_{k}\rangle
\end{equation}
of the residual vector $r_{k} = b - A x_{x}$ of the $k$-th iteration $x_{k}$. Algorithm 3 in \citet{jiang2020kaczmarz} and the Algorithm in \citet{steinerberger2020weighted} slightly differ with respect to the precise implementation, while Algorithm 2.1 in \cite{guan2020note} can be seen as a more restricted version. By referring to \citet{steinerberger2020weighted},
row $i$ is selected with probability
\begin{equation}\label{E4}
p_{k}(i) = \frac{|r_{k}(i)|^{p}}{\sum_{i=1}^{m} |r_{k}(i)|^{p}}\; .
\end{equation}
for some given integer $p$, see Algorithm \ref{A1}.

\begin{algorithm}
  \caption{Randomized Weighted Kaczmarz}\label{A1}
  Input: matrix $A\in \mathbb{C}^{m,n}$ with $m\geq n$ and $\|a_{i}\|_{2} = 1$ for all $i\in \{1,\ldots, m\}$, $b\in \mathbb{C}^{n}$ satisfying $b = A x$ for some $x$, initial guess $x_{0}$, integer $p$.
  \begin{itemize}
  \item[1.] Set $k=0$. Compute $r_{0} = b - A x_{0}$.
  \item[2.] Select an element $i$ from $\{1,\ldots, m\}$ with probability $p_{k}(i)$ from (\ref{E4}).
  \item[3.] Compute $x_{k+1} = x_{k} + r_{k}(i)\, a_{i}$ with $r_{k}(i)$ from (\ref{E3}).
  \item[4.] Compute $r_{k+1} = b - A x_{k}$. Set $k = k+1$. Go to 2.
  \end{itemize}
\end{algorithm}

As shown by \citet{steinerberger2020weighted}, Algorithm \ref{A1} is at least as efficient as the classical
randomized Kaczmarz method from \citet{strohmer2009randomized} irrespective of the choice of $p$, see also Algorithm 3 in \citet{jiang2020kaczmarz}. The case $p\rightarrow \infty$ corresponds to Algorithm \ref{A2} below, which actually is the maximum residual (MR) rule, coinciding with  the maximum distance (MD) rule for a standardized matrix $A$, see \cite{eldar2011acceleration, griebel2012greedy, nutini2016convergence}.

\begin{algorithm}
  \caption{Non-Randomized Greedy Kaczmarz}\label{A2}
  Input: matrix $A\in \mathbb{C}^{m,n}$ with $m\geq n$ and $\|a_{i}\|_{2} = 1$ for all $i\in \{1,\ldots, m\}$, $b\in \mathbb{C}^{n}$ satisfying $b = A x$ for some $x$, initial guess $x_{0}$.
  \begin{itemize}
  \item[1.] Set $k=0$. Compute $r_{0} = b - A x_{0}$.
  \item[2.] Select an element $i$ from $\{1,\ldots, m\}$ such that $\displaystyle i = \text{argmax}_{j} |r_{k}(j)|$ .
  \item[3.] Compute $x_{k+1} = x_{k} + r_{k}(i)\, a_{i}$.
  \item[4.] Compute $r_{k+1} = b - A x_{k}$. Set $k = k+1$. Go to 2.
  \end{itemize}
\end{algorithm}

The non-randomized Algorithm \ref{A2} is called `partially randomized' by \citet{jiang2020kaczmarz} and stated as their Algorithm 4. The authors show a convergence result that also relates Algorithm 4 to the greedy randomized Kaczmarz method
from \citet{bai2018greedy, bai2018relaxed}. It is noted, however, that there is nothing random about the actual selection of a row in Algorithm \ref{A2}, apart from the fact that a deterministic selection can be modelled as a probability one (almost sure) decision, see also the proof of the following Proposition in Appendix \ref{sec:A}.

\begin{prop} Let $x$ denote the solution of (\ref{E1}). Let $\mathbb{E}_{\text{Alg1}}$ and $\mathbb{E}_{\text{Alg2}}$ be the respective conditional expectation given $x_{k}$ corresponding to the probability distribution implied by the choice of row $i$ in step 2 of Algorithm \ref{A1} and \ref{A2}. Then
$$
E_{\text{Alg2}}\|x_{k+1} - x\|_{2}^{2} \leq E_{\text{Alg1}}\|x_{k+1} - x\|_{2}^{2}\; .
$$
\end{prop}

\begin{rmk}
The above Proposition together with the Theorem in \citet{steinerberger2020weighted} implies that Algorithm 2 converges with at least the rate of Algorithm 1.
\end{rmk}

\section{Partially Weighted Variant of Randomized Kaczmarz} \label{sec:greedy}

Both Algorithms \ref{A1} and \ref{A2} require the computation of the complete residual vector $r_{k}$ prior to the selection step 2 which can be quite time-consuming.
In order to reduce the number of required residuals $r_{k}(i)$, it has been proposed to select a random sample of rows in advance and then greedily select a row within the sample, see \cite{de2017sampling, haddock2021greed, jiang2020kaczmarz}.
Our Algorithm \ref{A3} below describes a proposal in the same light, where, however, the sample size by itself is random. It may also be seen as a variant of the classical randomized Kaczmarz method which makes use of individual residuals $r_{k}(i)$ only when needed, and does not require the computation of the complete residual vector $r_{k}$, except in rare cases.

The main idea is to look out for $\text{argmax}_{j} |r_{k}(j)|$, but with considerable reservation. In a first step a candidate row $i_{1}$ is randomly chosen from the set $\{1,\ldots , m\}$ according to the uniform distribution. Then a second competitor row $i_{2}$ is chosen from the remaining set $\{1,\ldots , m\}\setminus \{i_{1}\}$ due to the uniform distribution on this set. If $|r_{k}(i_{1})| > |r_{k}(i_{2})|$ row $i_{1}$ is selected, otherwise $i_{2}$ becomes the new candidate and a new competitor row $i_{3}$ is randomly selected from the set $\{1,\ldots , m\}\setminus \{i_{1},i_{2}\}$. The absolute residuals of the two rows in question are again compared in the same manner. This is repeated until a candidate row is actually selected, where one may always select the last possible row when all other residual comparisons did not lead to a candidate selection before.

The procedure requires at least the computation of 2 residuals for each iteration. At most, $m$ residuals are required, but this will rarely be the case. Of course, an obvious simpler variant is to compare only 2 residuals in each iteration and select the row admitting the larger residual. This may then be seen as employing a simple random sample of size 2. The performance of Algorithm \ref{A3} and its two residuals variant is considered in Section \ref{sec:num}.

Skipping steps 2.2 and 2.3 and employing $i = i_{1}$ from step 2.1 gives the classical randomized Kaczmarz by
\citet{strohmer2009randomized} for the case of a standardized matrix $A$.

\begin{rmk} If a specific row $i$ is selected using elements of $r_{k}$ by step 2 of one of the Algorithms \ref{A1}, \ref{A2}, and \ref{A3}, then $r_{k+1}(i) = 0$. Hence, for all three algorithms, a selected row in the subsequent iteration is (almost surely) different from the selected row in the actual iteration as long as the exact solution has not been found.
\end{rmk}

\begin{algorithm}
  \caption{Randomized Kaczmarz With Partially Weighted Selection Step}\label{A3}
  Input: matrix $A\in \mathbb{C}^{m,n}$ with $m\geq n$ and $\|a_{i}\|_{2} = 1$ for all $i\in \{1,\ldots, m\}$, $b\in \mathbb{C}^{n}$ satisfying $b = A x$ for some $x$, initial guess $x_{0}$.
  \begin{itemize}
  \item[1.] Set $k=0$.
  \item[2.] Set $U = \{1,\ldots, m\}$. Select an element $i$ from $U$ according to the following selection scheme:
  \begin{itemize}
  \item[2.1] Select an element $i_{1}$ from $U$ with uniform probability.
  \item[2.2] Set $U = U\setminus\{i_{1}\}$. If $U=\emptyset$ set $i = i_{1}$ and go to 3. Otherwise, select an element $i_{2}$ from the set $U$ with uniform probability.
  \item[2.3] If  $|r_{k}(i_{1})| > |r_{k}(i_{2})|$, set $i = i_{1}$ and go to 3. Otherwise, set $i_{1} = i_{2}$ and go to 2.2.
  \end{itemize}
  \item[3.] Compute $x_{k+1} = x_{k} + r_{k}(i)\, a_{i}$.
  \item[4.] Set $k = k+1$. Go to 2.
  \end{itemize}
\end{algorithm}

The convergence of Algorithm \ref{A3} can be concluded from the convergence of the classical randomized Kaczmarz algorithm. For each iteration the full set $\{1,\ldots, m\}$ of row numbers is again available for both, the classical and the partially weighted row selection step. Indeed, if there were a sequence of selected row numbers for which the partial weighted variant had failed to converge, then the classical algorithm too could not have converged for this very sequence.
See in addition \citet[Sect. 5.5]{patel2021convergence}, also confirming convergence of algorithms of this type.

Moreover, it is to be expected that Algorithm \ref{A3} converges faster than the classical randomized Kaczmarz, since the selection of rows is directed towards Algorithm \ref{A2}. Numerical examples discussed in the following section support this conclusion.

\section{Numerical Examples}\label{sec:num}

In this section we reconsider the settings described by \citet{steinerberger2020weighted}. Computations are carried out with the statistical software {\sf R}, see \citet{Rsoftware}.

\begin{figure}[hbt]
\centering
\includegraphics[width=19pc,angle=0]{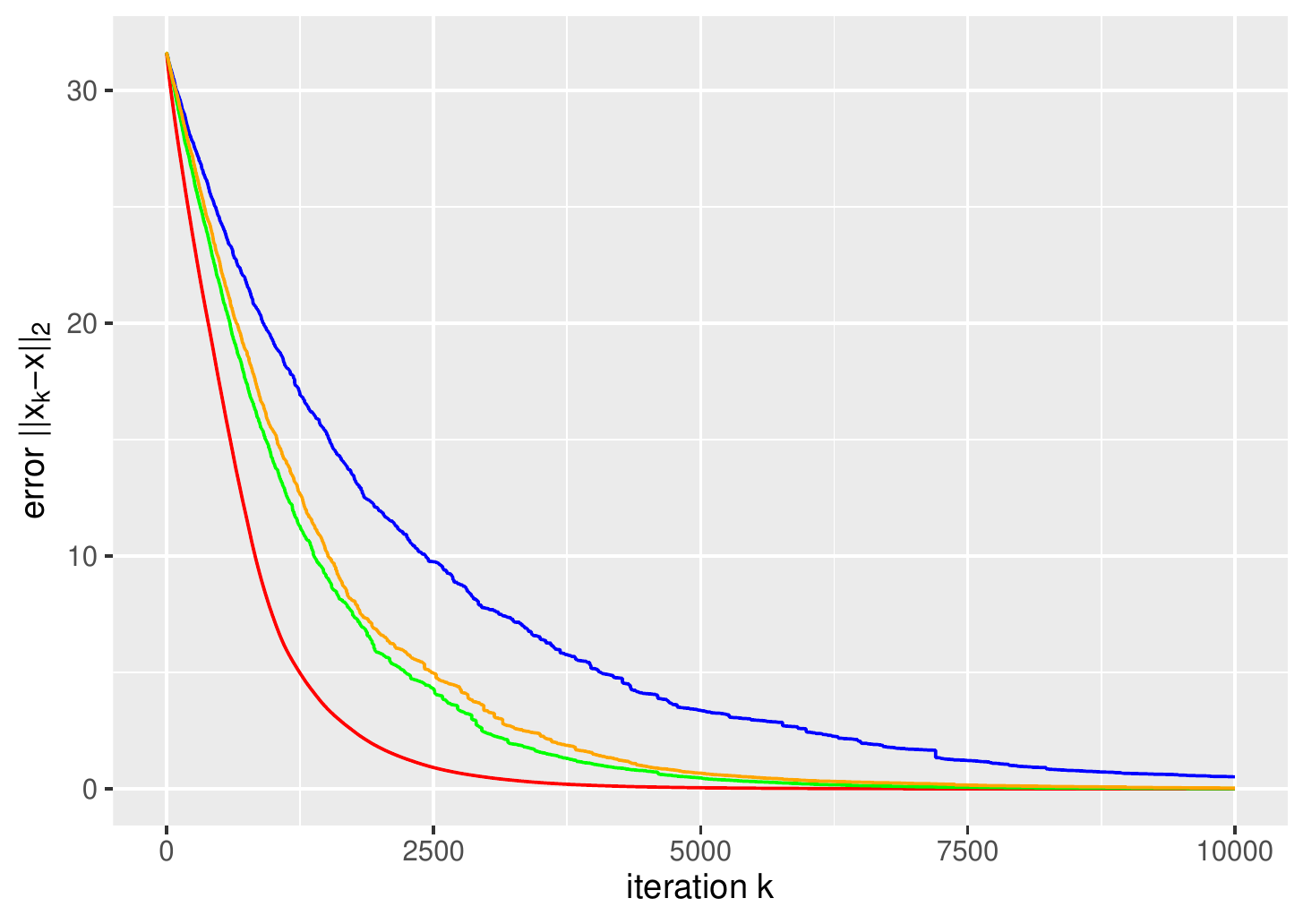}\\
  \caption{Error $\|x_{k} - x\|_{2}$ for the classical randomized Kaczmarz (blue), randomized Kaczmarz with partially weighted selection step (green), two residuals variant of the former (orange), and non-randomized greedy Kaczmarz (red)}\label{F1}
\end{figure}

\subsection{Nice Matrix} A $1000 \times 1000$ matrix $A$ is created by sampling the elements independently from the standard normal distribution. Then the matrix $100\, I_{1000}$ is added to $A$ and the result is standardized. The vector $b$ is the $1000\times 1$  vector of $0$s and the initial guess $x_{0}$ is the $1000\times 1$ vector of $1$s. The outcome from the considered algorithms is displayed in Figure \ref{F1}. The classical randomized Kaczmarz (blue line) performs in a very similar manner to what can be seen from Figure 2 in
\citet{steinerberger2020weighted}. Also, the non-randomized Kaczmarz from Algorithm \ref{A2} (red line) performs similar
to what can be seen from Figure 2 in \citet{steinerberger2020weighted} for the randomized weighted Kaczmarz (Algorithm \ref{A1}) for $p=20$. Visibly, the randomized Kaczmarz with partially weighted selection step (green line) performs better
than the classical randomized Kaczmarz and also slightly better than the two residuals variant. In order to assess the amount of complexity in Algorithm \ref{A3}, the number of required residuals is recorded in Table \ref{T1} for each of the first 10000 iterations.

\begin{table}[hbt]
\caption{Number of required residuals per iteration.}\label{T1}
\centering
\begin{tabular}{ccccccccc}
\hline
\# residuals & 2 &    3 &    4 &    5  &   6 &    7  &   8  &   9  \\
\hline
freq &4947 & 3334 & 1292  & 355   & 59   & 10   &  2  &    1\\
\hline
\end{tabular}
\end{table}

As it is seen, about half of the considered iterations require only the minimal number of two residuals, while no  iteration requires more than 9 residuals. Several repetitions have shown similar results.

\begin{figure}[hbt]
\centering
\includegraphics[width=19pc,angle=0]{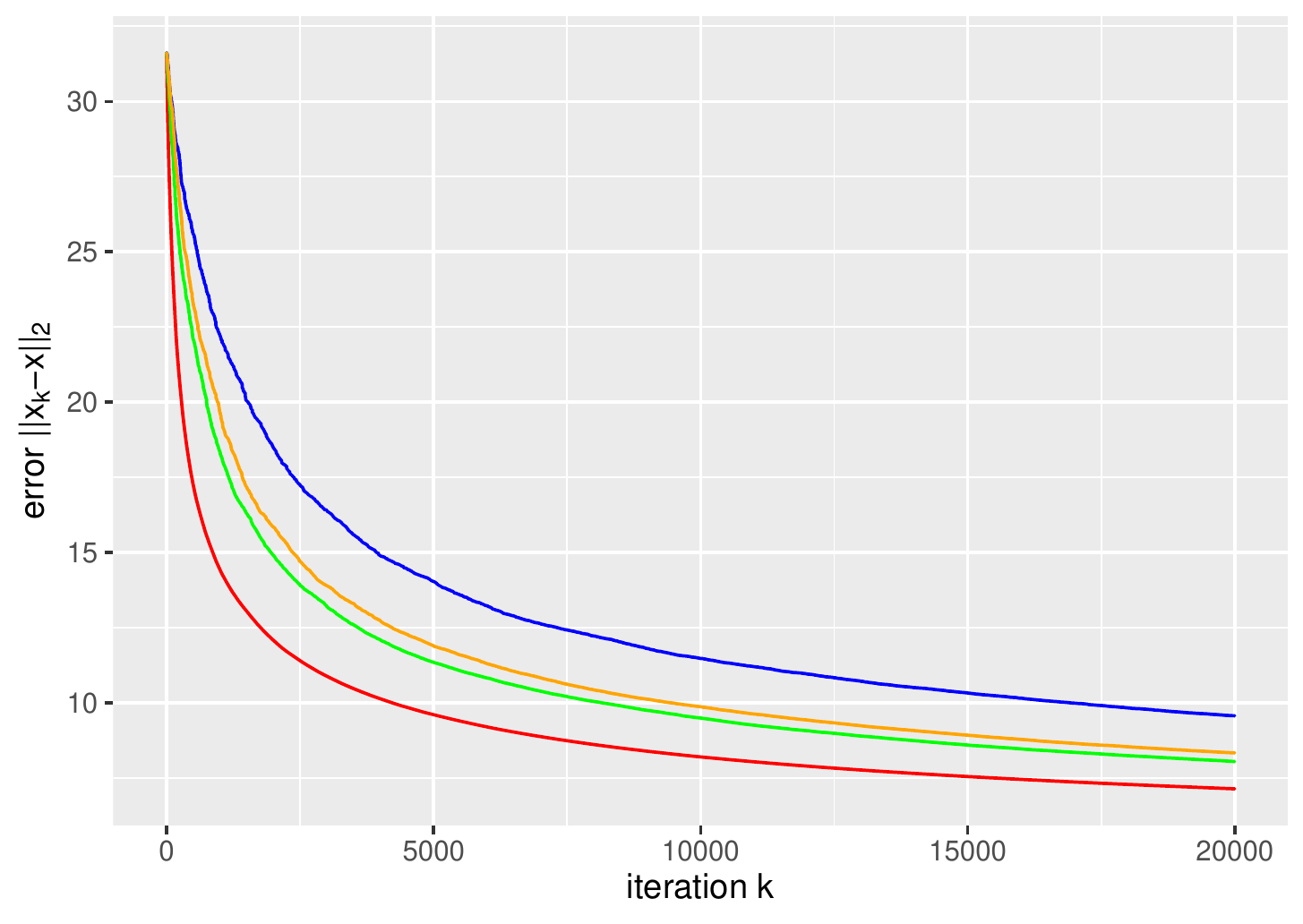}\\
  \caption{Error $\|x_{k} - x\|_{2}$ for the classical randomized Kaczmarz (blue), randomized Kaczmarz with partially weighted selection step (green), two residuals variant of the former (orange), and non-randomized greedy Kaczmarz (red)}\label{F2}
\end{figure}

\subsection{Challenging Matrix} For the more challenging setting, a $1000 \times 1000$ matrix $A$ is created by sampling the elements independently from the standard normal distribution. The vector $b$ and the initial choice $x_{0}$ are the same as before. Figure \ref{F2} displays the results. As before, results are comparable to those in
Figure 3 from \citet{steinerberger2020weighted}, the overall convergence being quite slower. Table \ref{T2} shows the number of required residuals in Algorithm  \ref{A3} for each of the first 20000 iterations. As in the case of the nice matrix, about half of the iterations require only two residuals and no iterations requires more than 9 residuals.

\begin{table}[hbt]
\caption{Number of required residuals per iteration.}\label{T2}
\centering
\begin{tabular}{ccccccccc}
\hline
\# residuals & 2 &    3 &    4 &    5  &   6 &    7  &   8  &   9  \\
\hline
freq & 9938 & 6665 & 2532  & 687  & 153 &   19 &    4 &    2\\
\hline
\end{tabular}
\end{table}

The results from the numerical examples support the conclusion that Algorithm \ref{A3} requires a rather moderate number of residuals, but has the capability to outperform the classical randomized Kaczmarz method.

\appendix

\section{Proof of the Proposition}\label{sec:A}

Let $a$ be a random vector with possible values $a_{1}, \ldots, a_{m}$, being the rows of the standardized matrix $A$ considered as column vectors. Let row $i$ be selected with $\text{prob}(i)$ for the computation of $x_{k+1}$. Let $d_{k} = x_{k} - x$, where $x$ denotes the solution to (\ref{E1}).

In the following, it is assumed that $x_{k}$ (and hence $d_{k}$) is given and thus non-random, while $x_{k+1}$ (and hence
$d_{k+1}$) is random, depending on the random vector $a$. The probability distribution $\text{prob}(i)$, $i=1,\ldots, m$, may be chosen as depending on $x_{k}$, and all considered expectations are regarded as conditional with respect to given $x_{k}$.

Now, the conditional expectation of the random variable $|\langle a, d_{k}\rangle|^2$ is given by
$$
\mathbb{E} |\langle a, d_{k}\rangle|^2 = \sum_{i=1}^{m} |\langle a_{i}, d_{k}\rangle|^2 \text{prob}(i)\; .
$$
Then, since the expectation of a discrete random variable cannot exceed its largest possible value,
$$
\mathbb{E} |\langle a, d_{k}\rangle|^2 \leq |\langle a_{i_{\ast}}, d_{k}\rangle|^2
$$
when $|\langle a_{i}, d_{k}\rangle|^2 \leq |\langle a_{i_{\ast}}, d_{k}\rangle|^2 $ for $i=1,\ldots ,m$.

Suppose that $\mathbb{E}_{\ast}$ denotes the conditional expectation of $|\langle a, d_{k}\rangle|^2$ with respect to the
specific probability distribution that assigns probability 1 to row $i_{\ast}$ satisfying
$$
|b_{i} - \langle a_{i}, x_{k}\rangle| \leq |b_{i_{\ast}} - \langle a_{i_{\ast}}, x_{k}\rangle|
$$
for $i=1,\ldots, m$. In view of the identity,
$$
|b_{i} - \langle a_{i}, x_{k}\rangle| = |\langle a_{i}, d_{k}\rangle|
$$
it follows that
$$
\mathbb{E}_{\ast}|\langle a, d_{k}\rangle|^2 = |\langle a_{i_{\ast}}, d_{k}\rangle|^2\; ,
$$
and thus
$$
\mathbb{E} |\langle a, d_{k}\rangle|^2 \leq \mathbb{E}_{\ast}|\langle a, d_{k}\rangle|^2\; .
$$
By following the lines in \citet[Sect. 4.1]{steinerberger2020weighted}, it is concluded that
$$
\|d_{k}\|^2 - |\langle a_{i}, d_{k}\rangle|^2, \quad i=1,\ldots , m\; ,
$$
are the possible values of the random variable $\|d_{k+1}\|_{2}^2$, so that
$$
\mathbb{E} \|d_{k+1}\|_{2}^2 = \|d_{k}\|_{2}^2 - \mathbb{E} |\langle a, d_{k}\rangle|^2
$$
is the conditional expectation of $\|d_{k+1}\|^2$ given $x_{k}$. Thus
$$
\mathbb{E}_{\ast} \|d_{k+1}\|_{2}^2 \leq \mathbb{E} \|d_{k+1}\|_{2}^2\; ,
$$
implying the Proposition.

\begin{rmk} The derived inequality $\mathbb{E}_{\ast} \|x_{k+1} - x\|_{2}^2 \leq \mathbb{E} \|x_{k+1} - x\|_{2}^2$ holds for any probability distribution specified for the selection of a row in computing $x_{k+1}$ from (\ref{E2}) when $A$ is standardized.
\end{rmk}


\begin{thebibliography}{29}
\providecommand{\natexlab}[1]{#1}
\providecommand{\url}[1]{\texttt{#1}}
\expandafter\ifx\csname urlstyle\endcsname\relax
  \providecommand{\doi}[1]{doi: #1}\else
  \providecommand{\doi}{doi: \begingroup \urlstyle{rm}\Url}\fi

\bibitem[Bai and Wu(2018{\natexlab{a}})]{bai2018greedy}
Z.-Z. Bai and W.-T. Wu.
\newblock On greedy randomized {K}aczmarz method for solving large sparse
  linear systems.
\newblock \emph{SIAM Journal on Scientific Computing}, 40:\penalty0 A592--A606,
  2018{\natexlab{a}}.

\bibitem[Bai and Wu(2018{\natexlab{b}})]{bai2018relaxed}
Z.-Z. Bai and W.-T. Wu.
\newblock On relaxed greedy randomized {K}aczmarz methods for solving large
  sparse linear systems.
\newblock \emph{Applied Mathematics Letters}, 83:\penalty0 21--26,
  2018{\natexlab{b}}.

\bibitem[De~Loera et~al.(2017)De~Loera, Haddock, and Needell]{de2017sampling}
J.~A. De~Loera, J.~Haddock, and D.~Needell.
\newblock A sampling {K}aczmarz--{M}otzkin algorithm for linear feasibility.
\newblock \emph{SIAM Journal on Scientific Computing}, 39:\penalty0 S66--S87,
  2017.

\bibitem[Deutsch(1984)]{deutsch1984rate}
F.~Deutsch.
\newblock Rate of convergence of the method of alternating projections.
\newblock In \emph{Parametric Optimization and Approximation}, pages 96--107.
  Springer, 1984.

\bibitem[Deutsch and Hundal(1997)]{deutsch1997rate}
F.~Deutsch and H.~Hundal.
\newblock The rate of convergence for the method of alternating projections,
  {II}.
\newblock \emph{Journal of Mathematical Analysis and Applications},
  205:\penalty0 381--405, 1997.

\bibitem[Eldar and Needell(2011)]{eldar2011acceleration}
Y.~C. Eldar and D.~Needell.
\newblock Acceleration of randomized {K}aczmarz method via the
  {J}ohnson--{L}indenstrauss lemma.
\newblock \emph{Numerical Algorithms}, 58:\penalty0 163--177, 2011.

\bibitem[Elfving et~al.(2014)Elfving, Hansen, and Nikazad]{elfving2014semi}
T.~Elfving, P.~C. Hansen, and T.~Nikazad.
\newblock Semi-convergence properties of {K}aczmarz’s method.
\newblock \emph{Inverse Problems}, 30\penalty0 (5):\penalty0 055007, 2014.

\bibitem[Gal\'{a}ntai(2005)]{galantai2005rate}
A.~Gal\'{a}ntai.
\newblock On the rate of convergence of the alternating projection method in
  finite dimensional spaces.
\newblock \emph{Journal of Mathematical Analysis and Applications},
  310\penalty0 (1):\penalty0 30--44, 2005.

\bibitem[Gordon et~al.(1970)Gordon, Bender, and Herman]{gordon1970algebraic}
R.~Gordon, R.~Bender, and G.~T. Herman.
\newblock Algebraic reconstruction techniques ({ART}) for three-dimensional
  electron microscopy and {X}-ray photography.
\newblock \emph{Journal of Theoretical Biology}, 29:\penalty0 471--481, 1970.

\bibitem[Gower et~al.(2019)Gower, Molitor, Moorman, and
  Needell]{gower2019adaptive}
R.~Gower, D.~Molitor, J.~Moorman, and D.~Needell.
\newblock Adaptive sketch-and-project methods for solving linear systems.
\newblock \emph{arXiv preprint arXiv:1909.03604}, 2019.

\bibitem[Gower and Richt\'{a}rik(2015)]{gower2015randomized}
R.~M. Gower and P.~Richt\'{a}rik.
\newblock Randomized iterative methods for linear systems.
\newblock \emph{SIAM Journal on Matrix Analysis and Applications}, 36:\penalty0
  1660--1690, 2015.

\bibitem[Griebel and Oswald(2012)]{griebel2012greedy}
M.~Griebel and P.~Oswald.
\newblock Greedy and randomized versions of the multiplicative schwarz method.
\newblock \emph{Linear Algebra and its Applications}, 437:\penalty0 1596--1610,
  2012.

\bibitem[Guan et~al.(2020)Guan, Li, Xing, and Qiao]{guan2020note}
Y.-J. Guan, W.-G. Li, L.-L. Xing, and T.-T. Qiao.
\newblock A note on convergence rate of randomized {K}aczmarz method.
\newblock \emph{Calcolo}, 57:\penalty0 1--11, 2020.

\bibitem[Haddock and Ma(2021)]{haddock2021greed}
J.~Haddock and A.~Ma.
\newblock Greed works: An improved analysis of sampling {K}aczmarz--{M}otzkin.
\newblock \emph{SIAM Journal on Mathematics of Data Science}, 3:\penalty0
  342--368, 2021.

\bibitem[Hefny et~al.(2017)Hefny, Needell, and Ramdas]{hefny2017rows}
A.~Hefny, D.~Needell, and A.~Ramdas.
\newblock Rows versus columns: Randomized {K}aczmarz or {G}auss--{S}eidel for
  ridge regression.
\newblock \emph{SIAM Journal on Scientific Computing}, 39\penalty0
  (5):\penalty0 S528--S542, 2017.

\bibitem[Jiang et~al.(2020)Jiang, Wu, and Jiang]{jiang2020kaczmarz}
Y.~Jiang, G.~Wu, and L.~Jiang.
\newblock A {K}aczmarz method with simple random sampling for solving large
  linear systems.
\newblock \emph{arXiv preprint arXiv:2011.14693}, 2020.

\bibitem[Kaczmarz(1937)]{karczmarz1937angenaherte}
S.~Kaczmarz.
\newblock {A}ngen\"{a}herte {A}ufl\"{o}sung von {S}ystemen linearer
  {G}leichungen.
\newblock \emph{Bull. Int. Acad. Polon. Sic. Lett.}, A:\penalty0 355--357,
  1937.

\bibitem[Liu and Wright(2016)]{liu2016accelerated}
J.~Liu and S.~Wright.
\newblock An accelerated randomized {K}aczmarz algorithm.
\newblock \emph{Mathematics of Computation}, 85:\penalty0 153--178, 2016.

\bibitem[Ma et~al.(2015)Ma, Needell, and Ramdas]{ma2015convergence}
A.~Ma, D.~Needell, and A.~Ramdas.
\newblock Convergence properties of the randomized extended {G}auss--{S}eidel
  and {K}aczmarz methods.
\newblock \emph{SIAM Journal on Matrix Analysis and Applications}, 36:\penalty0
  1590--1604, 2015.

\bibitem[Needell and Tropp(2014)]{needell2014paved}
D.~Needell and J.~A. Tropp.
\newblock Paved with good intentions: {A}nalysis of a randomized block
  {K}aczmarz method.
\newblock \emph{Linear Algebra and its Applications}, 441:\penalty0 199--221,
  2014.

\bibitem[Needell et~al.(2015)Needell, Zhao, and Zouzias]{needell2015randomized}
D.~Needell, R.~Zhao, and A.~Zouzias.
\newblock Randomized block {K}aczmarz method with projection for solving least
  squares.
\newblock \emph{Linear Algebra and its Applications}, 484:\penalty0 322--343,
  2015.

\bibitem[Nutini et~al.(2016)Nutini, Sepehry, Laradji, Schmidt, Koepke, and
  Virani]{nutini2016convergence}
J.~Nutini, B.~Sepehry, I.~Laradji, M.~Schmidt, H.~Koepke, and A.~Virani.
\newblock Convergence rates for greedy {K}aczmarz algorithms, and faster
  randomized {K}aczmarz rules using the orthogonality graph.
\newblock \emph{arXiv preprint arXiv:1612.07838}, 2016.

\bibitem[Patel et~al.(2021)Patel, Jahangoshahi, and
  Maldonado]{patel2021convergence}
V.~Patel, M.~Jahangoshahi, and D.~A. Maldonado.
\newblock Convergence of adaptive, randomized, iterative linear solvers.
\newblock \emph{arXiv preprint arXiv:2104.04816}, 2021.

\bibitem[Popa(2018)]{popa2018convergence}
C.~Popa.
\newblock Convergence rates for {K}aczmarz-type algorithms.
\newblock \emph{Numerical Algorithms}, 79:\penalty0 1--17, 2018.

\bibitem[{R Core Team}(2021)]{Rsoftware}
{R Core Team}.
\newblock \emph{R: A Language and Environment for Statistical Computing}.
\newblock R Foundation for Statistical Computing, Vienna, Austria, 2021.
\newblock URL \url{https://www.R-project.org/}.

\bibitem[Steinerberger(2020)]{steinerberger2020weighted}
S.~Steinerberger.
\newblock A weighted randomized {K}aczmarz method for solving linear systems.
\newblock \emph{arXiv preprint arXiv:2007.02910}, 2020.

\bibitem[Strohmer and Vershynin(2009)]{strohmer2009randomized}
T.~Strohmer and R.~Vershynin.
\newblock A randomized {K}aczmarz algorithm with exponential convergence.
\newblock \emph{Journal of Fourier Analysis and Applications}, 15:\penalty0
  262--278, 2009.

\bibitem[Zhang(2019)]{zhang2019new}
J.-J. Zhang.
\newblock A new greedy {K}aczmarz algorithm for the solution of very large
  linear systems.
\newblock \emph{Applied Mathematics Letters}, 91:\penalty0 207--212, 2019.

\bibitem[Zouzias and Freris(2013)]{zouzias2013randomized}
A.~Zouzias and N.~M. Freris.
\newblock Randomized extended {K}aczmarz for solving least squares.
\newblock \emph{SIAM Journal on Matrix Analysis and Applications}, 34:\penalty0
  773--793, 2013.

\end{thebibliography}
\end{document}